\documentclass[11pt]{article}
\usepackage{color,xcolor}
\usepackage{amsfonts} 
\usepackage{graphicx}
\usepackage{epstopdf}
\usepackage{amsmath}
\usepackage{amssymb}
\usepackage{multirow}
\usepackage{bbm}
\usepackage{color}
\usepackage{cite}
\usepackage{mathrsfs}

\newtheorem{theorem}{Theorem}[section]
\newtheorem{lemma}{Lemma}[section]

\newtheorem{definition}{Definition}[section]
\newtheorem{assumption}{Assumption}[section]
\newtheorem{example}{Example}[section]
\newtheorem{remark}{Remark}[section]

\numberwithin{equation}{section}

\begin{document}
\markboth{Jingna Zhang, Jianfei Huang, Yifa Tang, Luis V\'azquez}{A modified EM method and its fast implementation for multi-term Riemann-Liouville stochastic fractional differential equations}


\begin{center}{\bf A modified EM method and its fast implementation for multi-term Riemann-Liouville stochastic fractional differential equations}
\end{center}

\begin{center}
\centerline{\large Jingna Zhang$^{1,2}$, Jianfei Huang$^{3}$, Yifa Tang$^{1,2}$\footnote{Corresponding Author.~E-mail address:~tyf@lsec.cc.ac.cn(Y.F. Tang)}, Luis V\'azquez$^{4}$}
\vspace{0.25cm}
{\small {\em
$^1$LSEC, ICMSEC, Academy of Mathematics and Systems Science, Chinese Academy of Sciences, Beijing 100190, China\\
$^2$School of Mathematical Sciences, University of Chinese Academy of Sciences, Beijing 100049, China\\
$^3$College of Mathematical Sciences, Yangzhou University, Yangzhou 225002, China\\
$^4$Departamento de An\'alisis Matem\'atico y Matem\'atica Aplicada, Facultad de Inform\'atica, Universidad Complutense de Madrid, Madrid 28040, Spain}}
\end{center}
\begin{abstract}
In this paper, a modified Euler-Maruyama (EM) method is constructed for a kind of multi-term Riemann-Liouville stochastic fractional differential equations and the strong convergence order $\min\{1-\alpha_m, 0.5\}$ of the proposed method is proved with Riemann-Liouville fractional derivatives' orders $0<\alpha_1<\alpha_2<\cdots<\alpha_m<1$. Then, based on the sum-of-exponentials approximation, a fast implementation of the modified EM method which is called a fast EM method is derived to greatly improve the computational efficiency. Finally, some numerical examples are carried out to support the theoretical results and show the powerful computational performance of the fast EM method.
\end{abstract}

\begin{flushleft}
\end{flushleft}

\section{Introduction}	
As we all know, fractional differential equations (FDEs) can be viewed as a generalization of ordinary differential equations, which can be used to model complex physical phenomena and processes with nonlocal properties. In recent twenty years, with the continuous development of fractional differential equations, they have very important applications in the fields of mechanics, electrical engineering, electromagnetic wave, population system \cite{Sandev-2019,Tarasov-2009,Huang-2021,Rivero-2011,Luis-2018} and so on. Among them, multi-term FDEs are known as an important tool in describing viscoelastic damping materials, modelling nonlinear wave phenomenon in plasma and simulating anomalous diffusive process, such as anomalous relaxation in magnetic resonance imaging signal magnitude, mechanical models of oxygen delivery through capillaries \cite{Koeller-1984,Qin-2017,Srivastava-2010,Abdou-2017}, etc. Thus, there are lots of excellent researches focusing on the analytical methods and numerical methods for multi-term FDEs, see \cite{Kukla-2020,Daftardar-Gejji-2008,Dehghan-2015,Bu-2019,Gu-2017} and the references therein.

Meanwhile, with scientific research constantly deepening, researchers find that there are always some noise disturbances that could not be ignored in real life, and in order to better describe the phenomena and process influenced by these noisy factors, researchers pay more attention to stochastic differential equations (SDEs). As we all know, many classical SDEs play an extremely important role in describing physical phenomena. For instance, Langevin equations are used to describe Brownian motion and explain Einstein relations \cite{Kawasaki-1973} and the stochastic Navier-Stokes equation are often used to simulate various problems in fluid motion\cite{Temam-1979}. Nowadays, in addition to the field of physics, SDEs are widely used in option pricing, population growth and many other fields \cite{Deng-2020,Khodabin-2013,Podlubny-1999}. Especially, stochastic fractional differential equations (SFDEs) appeal more scholars' attention and many studies have been carried out, such as the random motion of harmonically trapped charged particles in a constant external magnetic field, the relationship between fluctuation-dissipation theorem and physical behavior, the noise driving in some financial models, epidemiological research \cite{Mankin-2018,Li-2017,Pedjeu-2012,Dung-2013} and so on. Obviously, it is quite difficult or even rarely impossible to obtain the exact solutions of SFDEs, thus there has been a growing interest to construct numerical methods for these equations. Until now, many numerical methods have been developed for SFDEs. For example, Kamrani in \cite{Kamrani-2015} investigated a numerical solution of SFDEs driven by additive noise and proved the convergence of the proposed method. Jin et al. in \cite{Jin-2019} studied the stochastic time-fractional diffusion problem driven by fractionally integrated Gaussian noise and developed a numerical scheme by employing the Galerkin finite element method in space, Gr$\ddot{u}$nwald-Letnikov formula in time and $L^2$-projection for the noise. Doan et al. in \cite{Doan-2020} constructed an Euler-Maruyama (EM) method for a kind of SFDEs driven by a multiplicative white noise with the Caputo fractional order $\alpha\in(0.5,1)$. Zhou et al. in \cite{Zhou-2021} used the finite difference method to solve the stochastic fractional nonlinear wave equation and presented the performance of numerical solution and property of energy under effect of two different types of noise, i.e., additive noise and multiplicative noise. Additionally, there also exist many other good works of SFDEs, see \cite{Anh-2019,Zheng-2020,Abouagwa-2019,Ahmadi-2017}, for examples.

Motivated by the above numerical methods for SFDEs, we in this paper will construct and analyze a modified EM method for the following multi-term Riemann-Liouville SFDEs,
\begin{eqnarray}\label{SFDE}
&&y'(t)+\sum_{i=1}^{m}\sideset{^{RL}_0}{^{\alpha_i}_t}{\mathop{\mathrm{D}}}y(t) = f(t,y(t)) + g(t,y(t))\frac{dW_t}{dt}, t \in [0,T],
\end{eqnarray}
where
\begin{itemize}
  \item $\sideset{^{RL}_0}{^{\alpha_i}_t}{\mathop{\mathrm{D}}}y(t)(i=1,2,\cdots,m)$ are the Riemann-Liouville fractional derivatives with $0<\alpha_1<\alpha_2<\cdots<\alpha_m<1$;
  \item $W_t$ is a one dimensional $\{\mathscr{F}_t\}_{\{0\leq t \leq T\}}$-adapted Brownian motion defined on the complete filtered probability space $\{\Omega,\mathscr{F},\mathbb{F}=\{\mathscr{F}_t\}_{t\in[0,\infty)},\mathbb{P}\}$;
  \item the initial value $y(0) = y_0\in \mathbb{R}^{d}$ is an $\mathscr{F}_0$-measurable random variable such that $E[|y_0|^2]<+\infty$ and $f,g:[0,T]\times\mathbb{R}^d\rightarrow\mathbb{R}^d$ are measurable functions.
\end{itemize}

In addition to the purpose of numerically solving the multi-term Riemann-Liouville SFDEs (\ref{SFDE}), we hope to construct more efficient numerical methods to avoid the huge computational cost caused by the nonlocal property of Riemann-Liouville fractional integral operators. Recently, lots of efficient methods based on the sum-of-exponentials (SOEs) technique proposed in \cite{Jiang-2017} by Jiang et al. are presented to overcome this difficulty for FDEs, see \cite{Zhu-2019,Cao-2020,Huang-2020,Sun-2021}, for examples. There are also some articles about applying the SOEs technique to efficiently solve the SFDEs. For example, Dai et al in \cite{Dai-2020} proposed an EM method and its fast implementation based on the SOEs approximation for L\'evy-driven stochastic Volterra integral equations with doubly singular kernels. Ma et al. \cite{Ma-2021} used the SOEs approximation to simulate the rough volatility and developed a fast two-step iteration algorithm. We in \cite{Zhang-2022} also presented a fast EM method for a class of nonlinear SFDEs by the SOEs technique. Herein, based on the SOEs approximation and the proposed EM method, a fast EM method for (\ref{SFDE}) is also introduced and analyzed in this paper.

The rest of this paper is organized as follows. In Section 2, some notations and preliminaries are given. In Section 3, the modified EM method is derived and the strong convergence of this method is proved. Section 4 aims to construct a fast EM method based on SOEs technique and present the corresponding numerical theoretical results. Two numerical examples are shown to examine the theoretical results and illustrate the effectiveness of the proposed two methods. Finally, a brief conclusion is given.

\section{Notations and preliminaries}
 In this section, we first introduce some useful notations that will be used throughout this paper. Let $E$ denote the expectation corresponding to $\mathbb{R}$. Define $\mathbb{L}^2(\Omega,\mathscr{F}_t,\mathbb{R})$ as the space of all $\mathscr{F}_t$-measurable, mean square integral functions $f=(f_1,f_2,\cdots,f_d)^{\top}:\Omega\rightarrow\mathbb{R}^{d}$ with standard Euclidean norm $\|f\|=\sqrt{\sum_{i=1}^{d}E\left[|f_{i}|^2\right]}$. For two real numbers $A$ and $B$, we denote $\max\{A,B\} = A \vee B$ and $\min\{A,B\} = A \wedge B$. Moreover, the capital letter $C$ will be used to represent a positive constant whose value may change when it appears in different places. The following two definitions are taken from \cite{Mao-2007}.
\begin{definition}\label{Riemann-Liouville}
The Riemann-Liouville fractional integral of order $\alpha (\alpha \geq 0)$ of function $f:[0,+\infty)\rightarrow \mathbb{R}^{d}$ is defined by
\begin{small}
\begin{eqnarray*}
\sideset{^{}_0}{^{\alpha}_t}{\mathop{\mathrm{J}}}f(t) = \frac{1}{\Gamma(\alpha)}\int_{0}^{t}(t-\tau)^{\alpha-1}f(\tau)d\tau,
\end{eqnarray*}
\end{small}
with ${\mathop{\mathrm{J}}}^{0}f(t) = f(t)$, where $\Gamma(\alpha)=\int_{0}^{+\infty}e^{-t}t^{\alpha-1}dt$ is the Euler gamma function.
\end{definition}

\begin{definition}\label{Riemann-Liouville_derivative}
The Riemann-Liouville derivative of order $\alpha(0\leq\alpha<1)$ of function $f \in C([0,+\infty))$ can be written as
\begin{small}
\begin{eqnarray*}
\sideset{^{RL}_0}{^{\alpha}_t}{\mathop{\mathrm{D}}} f(t) = \frac{d}{dt}\left(\frac{1}{\Gamma(1-\alpha)}\int_{0}^{t}\frac{f(\tau)}{(t-\tau)^{\alpha}}d\tau\right).
\end{eqnarray*}
\end{small}
\end{definition}
In the rest of this section, to guarantee the existence and uniqueness of the solution of equation (\ref{SFDE}), we impose three important hypotheses.
\begin{assumption}\label{assumption_1}
(H$\ddot{o}$lder continuity) There exists a positive constant $L_1$ such that for all $t_1,t_2\in[0,T]$ and $y\in\mathbb{R}^d$, $f$ and $g$ satisfy the condition:
\begin{eqnarray*}
\|f(t_1,y)-f(t_2,y)\| \vee \|g(t_1,y)-g(t_2,y)\|\leq L_1|t_1-t_2|.
\end{eqnarray*}
\end{assumption}

\begin{assumption}\label{assumption_2}
(Lipschitz continuity) There exists a positive constant $L_2$ such that for all $t\in[0,T]$ and $x_1,x_2\in\mathbb{R}^d$, $f$ and $g$ satisfy the inequality:
\begin{eqnarray*}
\|f(t,x_1)-f(t,x_2)\|\vee \|g(t,x_1)-g(t,x_2)\| \leq L_2\|x_1-x_2\|.
\end{eqnarray*}
\end{assumption}

\begin{assumption}\label{assumption_3}
(Linear growth condition) There exists a positive constant $L_3$ such that for all $t\in[0,T]$ and all $y\in\mathbb{R}^d$, $f$ and $g$ satisfy the condition:
\begin{eqnarray*}
\|f(t,y)\|\vee\|g(t,y)\| \leq L_3(1+\|y\|).
\end{eqnarray*}
\end{assumption}

\section{The modified EM method for multi-term Riemann-Liouville SFDEs}
\subsection{The modified EM method}
In order to get the numerical method of the Riemann-Liouville SFDEs (\ref{SFDE}), we first transform equation (\ref{SFDE}) into the following stochastic Volterra integral equation (SVIE) by taking the Riemann-Liouville fractional integral operator $\sideset{^{}_0}{^{}_t}{\mathop{\mathrm{J}}}$ on both sides of equation (\ref{SFDE}), that is
\begin{small}
\begin{eqnarray}\label{SVIE}
y(t) &=& y_0-\sum_{i=1}^{m}\frac{1}{\Gamma(1-\alpha_i)}\int_{0}^{t}(t-s)^{-\alpha_i}y(s)ds + \int_{0}^{t}f(s,y(s))ds + \int_{0}^{t}g(s,y(s))dW_s,
\end{eqnarray}
\end{small}
where the property of Riemann-Liouville fractional integral and Riemann-Liouville fractional derivative is used, i.e.
\begin{small}
\begin{eqnarray*}
\sideset{^{}_0}{^{}_t}{\mathop{\mathrm{J}}}(\sideset{^{RL}_0}{^{\alpha}_t}{\mathop{\mathrm{D}}}y(t)) = \frac{1}{\Gamma(1-\alpha_i)}\int_{0}^{t}\frac{y(s)}{(t-s)^{\alpha_i}}ds.
\end{eqnarray*}
\end{small}
In another way, for each $y_0 \in \mathbb{L}^2(\Omega,\mathscr{F}_t,\mathbb{R})$, when the above equality holds for $t\in[0,T]$, a $\mathbb{F}$-adapted process $y(t)$ is called a solution of equation (\ref{SFDE}) on the interval $[0,T]$ with the initial condition $y(0)=y_0$.

For every integer $N \geq 1$, the modified EM method for equation (\ref{SVIE}) can be presented as
\begin{small}
\begin{eqnarray}\label{the discrete EM scheme}
Y^{(N)}(t_n) &=& y_0 - \sum_{i=1}^{m}\frac{1}{\Gamma(1-\alpha_i)}\sum_{j=0}^{n-1}(t_n-t_j)^{-\alpha_i}Y^{(N)}(t_j)h \nonumber\\ &&+ \sum_{j=0}^{n-1}f(t_j,Y^{(N)}(t_j))h + \sum_{j=0}^{n-1}g(t_j,Y^{(N)}(t_j))\triangle W_j,
\end{eqnarray}
\end{small}
where $t_n = nh (n = 0, 1, 2, \cdots, N)$ denote the grid points with step size $h = \frac{T}{N}$, $\triangle W_j = W(t_{j+1})-W(t_j) (j = 0, 1, 2, \cdots, N-1)$ denote the increments of Brownian motion. To facilitate its convergence analysis, we introduce a continuous-time version
\begin{small}
\begin{eqnarray}\label{EM-scheme}
Y^{(N)} (t) &=& y_0 - \sum_{i=1}^{m}\frac{1}{\Gamma(1-\alpha_i)}\int_{0}^{t}(t-\tau_N(s))^{-\alpha_i}Y^{(N)}(\tau_N(s))ds\nonumber\\
&&+\int_{0}^{t}f(\tau_N(s),Y^{(N)}(\tau_N(s)))ds +\int_{0}^{t}g(\tau_N(s),Y^{(N)}(\tau_N(s)))dW(s),
\end{eqnarray}
\end{small}
where $\tau_N(s) = \frac{nT}{N}, n = 0, 1, \cdots, N-1$ for $s \in (\frac{nT}{N},\frac{(n+1)T}{N}]$.

\subsection{Strong convergence of the modified EM method}
To prove the strong convergence of the modified EM method (\ref{EM-scheme}), we first list some necessary lemmas.
\begin{lemma}\label{lemma4.1}
 If $\alpha\in(0,1]$, then for any $s' \leq s < \tilde{t} < t$, the following inequality holds
 \begin{small}
\begin{eqnarray*}
\frac{1}{(\tilde{t}-s')^{\alpha}}-\frac{1}{(t-s')^{\alpha}} \leq \frac{1}{(\tilde{t}-s)^{\alpha}}-\frac{1}{(t-s)^{\alpha}}.
\end{eqnarray*}
\end{small}
\end{lemma}
\noindent\emph{Proof.}
Define $p(x) = \frac{1}{(\tilde{t}-x)^{\alpha}}-\frac{1}{(t-x)^{\alpha}}$ with $x < \tilde{t} < t$. Then we have $p'(x) = \alpha\left[\frac{1}{(\tilde{t}-s)^{\alpha+1}}-\frac{1}{(t-s)^{\alpha+1}}\right]$. According to the property of power function, it holds that $p'(x) > 0$. Thus it implies that $p(s') \leq p(s)$ together with $s' \leq s$.
\ \ $\Box$

\begin{lemma}\label{lemma4.2}
 Suppose $\alpha\in(0,1)$, then there exists a positive constant $C$ independent of $h$ such that for any $t\in(t_n,t_{n+1}], n=0,1,2,\cdots, N-1$,
 \begin{small}
\begin{eqnarray*}
\int_{0}^{t}|(t-\tau_N(s))^{-\alpha}-(t-s)^{-\alpha}|ds \leq Ch^{1-\alpha},
\end{eqnarray*}
\end{small}
where $\tau_N(s)$ have been defined in the above subsection.
\end{lemma}
\noindent\emph{Proof.}
For $t\in(t_n,t_{n+1}]$, it is easily derived that
\begin{small}
\begin{eqnarray*}
&&\int_{0}^{t}|(t-\tau_N(s))^{-\alpha}-(t-s)^{-\alpha}|ds\\
&=& \int_{0}^{t}(t-s)^{-\alpha}-(t-\tau_N(s))^{-\alpha}ds\\
&\leq& \sum_{j=0}^{n-1}\int_{t_j}^{t_{j+1}}(t_n-s)^{-\alpha}-(t_{n+1}-t_j)^{-\alpha}ds + \int_{t_n}^{t}(t-s)^{-\alpha}-(t-t_n)^{-\alpha}ds\\
&\leq&\frac{t_n}{1-\alpha}-\sum_{j=0}^{n-1}(t_{n+1}-t_j)^{-\alpha}h + \frac{\alpha}{1-\alpha}(t-t_n)^{1-\alpha}\\
&\leq&\frac{(nh)^{1-\alpha}}{1-\alpha} - n((n+1)h)^{-\alpha}h + \frac{\alpha}{1-\alpha}h^{1-\alpha}\\
&\leq&Ch^{1-\alpha}.
\end{eqnarray*}
\end{small}
\ \ $\Box$

The following inequality can easily be obtained by the property of norm.
\begin{lemma}\label{inequality}
For all $y_i\in\mathbb{R}^d$, $i=1,2,\cdots,n$, the following inequality holds
\begin{small}
\begin{eqnarray*}
\|\sum_{i=1}^{n}y_i\|^{2} \leq n\sum_{i=1}^{n}\|y_{i}\|^{2}.
\end{eqnarray*}
\end{small}
\end{lemma}

 To carry out the proof of the strong convergence of the modified EM method (\ref{EM-scheme}), the bounded estimate for the numerical solutions is given.

 \begin{theorem}\label{bound}
 Suppose Assumption \ref{assumption_3} holds, then for all $N\in N^{*}$, the numerical solution $Y^{(N)}(t)$ of the modified EM method (\ref{EM-scheme}) is bounded, that is
 \begin{footnotesize}
 \begin{eqnarray*}
E\left[\|Y^{(N)}(t)\|^2\right] \leq C \quad {\rm{and}} \quad E\left[\|Y^{(N)}(t_n)\|^2\right] \leq C, \quad \forall t, t_n \in [0,T],
\end{eqnarray*}
\end{footnotesize}
where $C$ is a positive constant independent of $N$.
 \end{theorem}

\noindent\emph{Proof.}
From equation (\ref{EM-scheme}) and Lemma \ref{inequality}, it is deduced that
\begin{small}
\begin{eqnarray*}\label{bounded_1}
&&\frac{1}{4}E\left[\|Y^{(N)}(t)\|^2\right]\\
 &\leq& E\left[\|y_0\|^2\right] + E\left[\|\sum_{i=1}^{m}\frac{1}{\Gamma(1-\alpha_i)}\int_{0}^{t}(t-\tau_N(s))^{-\alpha_i}Y^{(N)}(\tau_N(s))ds\|^2\right]\\
&&+E\left[\|\int_{0}^{t}f(\tau_N(s),Y^{(N)}(\tau_N(s)))ds\|^2\right] + E\left[\|\int_{0}^{t}g(\tau_N(s),Y^{(N)}(\tau_N(s)))dW(s)\|^2\right].
\end{eqnarray*}
\end{small}
Then using the H$\ddot{o}$lder inequality and It$\hat{o}$ isometry as well as Lemma \ref{inequality} implies that
\begin{small}
\begin{eqnarray*}
&&\frac{1}{4}E\left[\|Y^{(N)}(t)\|^2\right]\\
 &\leq& E\left[\|y_0\|^2\right] + \sum_{i=1}^{m}\frac{m}{\Gamma^2(1-\alpha_i)}\int_{0}^{t}\frac{1}{(t-\tau_N(s))^{\alpha_i}}ds\int_{0}^{t}\frac{E\left[\|Y^{(N)}(\tau_N(s))\|^2\right]}{(t-\tau_N(s))^{\alpha_i}}ds\\
&&+t\int_{0}^{t}E\left[\|f(\tau_N(s),Y^{(N)}(\tau_N(s)))\|^2\right]ds + \int_{0}^{t}E\left[\|g(\tau_N(s),Y^{(N)}(\tau_N(s)))\|^2\right]ds.
\end{eqnarray*}
\end{small}
 Applying the linear growth condition of Assumption \ref{assumption_3} and the property of power function with the fact $\tau_N(s) \leq s$, it can be deduced that
\begin{small}
\begin{eqnarray*}
&&\frac{1}{4}E\left[\|Y^{(N)}(t)\|^2\right]\\
 &\leq& E\left[\|y_0\|^2\right] + \sum_{i=1}^{m}\frac{m}{\Gamma^2(1-\alpha_i)}\int_{0}^{t}\frac{1}{(t-s)^{\alpha_i}}ds\int_{0}^{t}\frac{E\left[\|Y^{(N)}(\tau_N(s))\|^2\right]}{(t-s)^{\alpha_i}}ds\\
&&+t\int_{0}^{t}2L_3^2(1+E\left[\|Y^{(N)}(\tau_N(s))\|^2\right])ds + \int_{0}^{t}2L_3^2(1+E\left[\|Y^{(N)}(\tau_N(s))\|^2\right])ds.
\end{eqnarray*}
\end{small}
Arranging the above inequality, it becomes
\begin{small}
\begin{eqnarray*}
&&E\left[\|Y^{(N)}(t)\|^2\right]\\
 &\leq& 4E\left[\|y_0\|^2\right] + \sum_{i=1}^{m}\frac{4m}{\Gamma^2(1-\alpha_i)}\frac{t^{1-\alpha_i}}{1-\alpha_i}\int_{0}^{t}(t-s)^{-\alpha_i}E\left[\|Y^{(N)}(\tau_N(s))\|^2\right]ds\\
&&+(8t+8)L_3^2\int_{0}^{t}1+E\left[\|Y^{(N)}(\tau_N(s))\|^2\right]ds\\
&\leq&4E\left[\|y_0\|^2\right] + (8T^2+8T)L_3^2 + (8T+8)L_3^2\int_{0}^{t}E\left[\|Y^{(N)}(\tau_N(s))\|^2\right]ds\\
&&+ \sum_{i=1}^{m}\frac{4mT^{1-\alpha_i}}{\Gamma(1-\alpha_i)\Gamma(2-\alpha_i)}\int_{0}^{t}(t-s)^{-\alpha_i}E\left[\|Y^{(N)}(\tau_N(s))\|^2\right]ds.
\end{eqnarray*}
\end{small}
There exists a positive constant $\alpha\in\{\alpha_1,\alpha_2,\cdots,\alpha_m\}$ such that
\begin{small}
\begin{eqnarray*}
E\left[\|Y^{(N)}(t)\|^2\right] &\leq& 4E\left[\|y_0\|^2\right] + (8T^2+8T)L_3^2 + C\int_{0}^{t}(t-s)^{-\alpha}E\left[\|Y^{(N)}(\tau_N(s))\|^2\right]ds
\end{eqnarray*}
\end{small}
Taking the supremum on both sides of the above inequality and according to the Gronwall's inequality, we can arrive at
\begin{small}
\begin{eqnarray*}
E\left[\|Y^{(N)}(t)\|^2\right] \leq C.
\end{eqnarray*}
\end{small}
Similarly, we can prove that $E\left[\|Y^{(N)}(t_n)\|^2\right] \leq C$. This proof is completed.
\ \ $\Box$

\begin{lemma}\label{t_t_n}
Suppose Assumption \ref{assumption_3} holds, then for all $N\in N^{*}$, there exists a positive constant $C$ independent of $N$ such that for all $t\in(t_n,t_{n+1}], t_n\in[0,T]$,
\begin{small}
\begin{eqnarray*}
E\left[\|Y^{(N)}(t)-Y^{(N)}(t_n)\|^2\right] \leq Ch^{{\min\{2(1-\alpha_m),1\}}}.
\end{eqnarray*}
\end{small}
\end{lemma}

\noindent\emph{Proof.}
For arbitrary $t\in(t_n,t_{n+1}], t_n\in[0,T]$, it follows from equation (\ref{EM-scheme}) that
\begin{small}
\begin{eqnarray*}
&&Y^{(N)}(t)-Y^{(N)}(t_n)\\
 &=& -\sum_{i=1}^{m}\frac{1}{\Gamma(1-\alpha_i)}\left\{\int_{0}^{t_n}\left[(t-\tau_N(s))^{-\alpha_i}-(t_n-\tau_N(s))^{-\alpha_i}\right]Y^{(N)}(\tau_N(s))ds\right.\\
 &&\left.+\int_{t_n}^{t}(t-\tau_N(s))^{-\alpha_i}Y^{(N)}(\tau_N(s))ds\right\}\\
&&+\int_{t_n}^{t}f(\tau_N(s),Y^{(N)}(\tau_N(s)))ds + \int_{t_n}^{t}g(\tau_N(s),Y^{(N)}(\tau_N(s)))dW(s).
\end{eqnarray*}
\end{small}
Applying Lemma \ref{inequality}, we obtain
\begin{small}
\begin{eqnarray*}
&&\frac{1}{4}E\left[\|Y^{(N)}(t)-Y^{(N)}(t_n)\|^2\right]\\
 &\leq& m\sum_{i=1}^{m}E\left[\|\frac{1}{\Gamma(1-\alpha_i)}\int_{0}^{t_n}\left[(t-\tau_N(s))^{-\alpha_i}-(t_n-\tau_N(s))^{-\alpha_i}\right]Y^{(N)}(\tau_N(s))ds\|^2\right]\\
 &&+m\sum_{i=1}^{m}E\left[\|\frac{1}{\Gamma(1-\alpha_i)}\int_{t_n}^{t}(t-\tau_N(s))^{-\alpha_i}Y^{(N)}(\tau_N(s))ds\|^2\right]\\
 &&+E\left[\|\int_{t_n}^{t}f(\tau_N(s),Y^{(N)}(\tau_N(s)))ds\|^2\right]+ E\left[\|g(\tau_N(s),Y^{(N)}(\tau_N(s)))dW(s)\|^2\right].
\end{eqnarray*}
\end{small}
By H$\ddot{o}$lder inequality and It$\hat{o}$ isometry, the above inequality turns into
\begin{small}
\begin{eqnarray*}
&&\frac{1}{4}E\left[\|Y^{(N)}(t)-Y^{(N)}(t_n)\|^2\right]\\
 &\leq& \sum_{i=1}^{m}\frac{m}{\Gamma^2(1-\alpha_i)}\int_{0}^{t_n}|(t-\tau_N(s))^{-\alpha_i}-(t_N-\tau_N(s))^{-\alpha_i}|ds\\
 &&\cdot\int_{0}^{t_n}|(t-\tau_N(s))^{-\alpha_i}-(t_n-\tau_N(s))^{-\alpha_i}|E\left[\|Y^{(N)}(\tau_N(s))\|^2\right]ds\\
 &&+\sum_{i=1}^{m}\frac{m}{\Gamma^2(1-\alpha_i)}\int_{t_n}^{t}(t-\tau_N(s))^{-\alpha_i}ds\int_{t_n}^{t}(t-\tau_N(s))^{-\alpha_i}E\left[\|Y^{(N)}(\tau_N(s))\|^2\right]ds\\
 &&+(t-t_n)\int_{t_n}^{t}E\left[\|f(\tau_N(s),Y^{(N)}(\tau_N(s)))\|^2\right]ds + \int_{t_n}^{t}E\left[\|g(\tau_N(s),Y^{(N)}(\tau_N(s)))\|^2\right]ds.
\end{eqnarray*}
\end{small}
This together with Lemma \ref{lemma4.1}, Assumption \ref{assumption_3} and Theorem \ref{bound} as well as the fact $\tau_N(s) \leq s$ implies that
\begin{small}
\begin{eqnarray*}
&&\frac{1}{4}E\left[\|Y^{(N)}(t)-Y^{(N)}(t_n)\|^2\right]\\
 &\leq& \sum_{i=1}^{m}\frac{Cm}{\Gamma^2(1-\alpha_i)}\int_{0}^{t_n}(t_n-s)^{-\alpha_i}-(t-s)^{-\alpha_i}ds\int_{0}^{t_n}(t_n-s)^{-\alpha_i}-(t-s)^{-\alpha_i}ds\\
 &&+\sum_{i=1}^{m}\frac{Cm}{\Gamma^2(1-\alpha_i)}\int_{t_n}^{t}(t-s)^{-\alpha_i}ds\int_{t_n}^{t}(t-s)^{-\alpha_i}ds + 2L_3^2(1+C)\left[(t-t_n)^{2} + (t-t_n)\right] \\
 &\leq&\sum_{i=1}^{m}\frac{Cm}{\Gamma^2(2-\alpha_i)}\left[(t-t_n)^{1-\alpha_i}+t_n^{1-\alpha_i}-t^{1-\alpha_i}\right]^2 + \sum_{i=1}^{m}\frac{mC}{\Gamma^2(2-\alpha_i)}(t-t_n)^{2-2\alpha_i}\\
 &&+2L_3^2(1+C)\left[(t-t_n)^2+(t-t_n)\right].
\end{eqnarray*}
\end{small}
Noticing $t-t_n \leq h$, thus we can get
\begin{eqnarray*}
E\left[\|Y^{(N)}(t)-Y^{(N)}(t_n)\|^2\right] \leq Ch^{{\min\{2(1-\alpha_m),1\}}}.
\end{eqnarray*}
This proof is completed.
\ \ $\Box$

\begin{theorem}\label{convergence}(Strong convergence)
Under Assumptions \ref{assumption_1}-\ref{assumption_3}, there exists a positive constant $C$ independent of $N$ such that for all $t\in[0,T]$,
\begin{small}
\begin{eqnarray*}
E\left[\|Y^{(N)}(t)-y(t)\|^2\right] \leq Ch^{{\min\{2(1-\alpha_m), 1\}}}.
\end{eqnarray*}
\end{small}
\end{theorem}

\noindent\emph{Proof.}
From SVIE (\ref{SVIE}) and the modified EM method (\ref{EM-scheme}), by adding some intermediate items and using Lemma \ref{inequality}, H$\ddot{o}$lder inequality and It$\hat{o}$ isometry, we derive that
\begin{footnotesize}
\begin{eqnarray*}
&&\frac{1}{3}E\left[\|Y^{(N)}(t)-y(t)\|^2\right] \\
&\leq&\sum_{i=1}^{m}\frac{3m}{\Gamma^2(1-\alpha_i)}\left\{\int_{0}^{t}|(t-\tau_N(s))^{-\alpha_i}-(t-s)^{-\alpha_i}|ds\right.\\
&&\left.\cdot\int_{0}^{t}|(t-\tau_N(s))^{-\alpha_i}-(t-s)^{-\alpha_i}|E\left[\|Y^{(N)}(\tau_N(s))\|^2\right]ds\right.\\
&&\left.+\int_{0}^{t}(t-s)^{-\alpha_i}ds\int_{0}^{t}(t-s)^{-\alpha_i}E\left[\|Y^{(N)}(\tau_N(s))-Y^{(N)}(s)\|^2\right]ds\right.\\
&&\left.+\int_{0}^{t}(t-s)^{-\alpha_i}ds\int_{0}^{t}(t-s)^{-\alpha_i}E\left[\|Y^{(N)}(s)-y(s)\|^2\right]ds\right\}\\
&&+3t\left\{\int_{0}^{t}E\left[\|f(\tau_N(s),Y^{(N)}(\tau_N(s)))-f(s,Y^{(N)}(\tau_N(s)))\|^2\right]ds\right.\\
&&\left.+\int_{0}^{t}E\left[\|f(s,Y^{(N)}(\tau_N(s)))-f(s,Y^{(N)}(s))\|^2\right]ds + \int_{0}^{t}E\left[\|f(s,Y^{(N)}(s))-f(s,y(s))\|^2\right]ds\right\}\\
&&+3\left\{\int_{0}^{t}E\left[\|g(\tau_N(s),Y^{(N)}(\tau_N(s)))-g(s,Y^{(N)}(\tau_N(s)))\|^2\right]ds\right.\\
&&\left.+\int_{0}^{t}E\left[\|g(s,Y^{(N)}(\tau_N(s)))-g(s,Y^{(N)}(s))\|^2\right]ds + \int_{0}^{t}E\left[\|g(s,Y^{(N)}(s))-g(s,y(s))\|^2\right]ds\right\}.
\end{eqnarray*}
\end{footnotesize}
Based on Assumptions \ref{assumption_1} and \ref{assumption_2}, Lemmas \ref{lemma4.2} and \ref{t_t_n} and Theorem \ref{bound} as well as the fact $s-\tau_N(s) \leq h$, it is easy to obtain that
\begin{small}
\begin{eqnarray*}
&&\frac{1}{3}E\left[\|Y^{(N)}(t)-y(t)\|^2\right] \\
&\leq&\sum_{i=1}^{m}\frac{3mC}{\Gamma^2(1-\alpha_i)}\left\{h^{2(2-\alpha_i)} + h^{\min\{2(1-\alpha_m),1\}} + \int_{0}^{t}\frac{E\left[\|Y^{(N)}(s)-y(s)\|^2\right]}{(t-s)^{\alpha_i}}ds\right\}\\
&&+3t\left\{L_1^2h^2t + L_2^2Ch^{{\min\{2(1-\alpha_m),1\}}} + L_2^2\int_{0}^{t}E\left[\|Y^{(N)}(s)-y(s)\|^2\right]ds\right\}\\
&&+3\left\{L_1^2h^2t + L_2^2Ch^{{\min\{2(1-\alpha_m),1\}}} + L_2^2\int_{0}^{t}E\left[\|Y^{(N)}(s)-y(s)\|^2\right]ds\right\}.
\end{eqnarray*}
\end{small}
Arranging the above inequality yields
\begin{small}
\begin{eqnarray*}
E\left[\|Y^{(N)}(t)-y(t)\|^2\right] \leq Ch^{{\min\{2(1-\alpha_m,1)\}}} + C\int_{0}^{t}(t-s)^{-\alpha}E\left[\|Y^{(N)}(s)-y(s)\|^2\right]ds,
\end{eqnarray*}
\end{small}
where $\alpha\in\{\alpha_1,\alpha_2,\cdots,\alpha_m\}$. Then using the Gronwall's inequality and arbitrariness of $t\in[0,T]$, we arrive at
\begin{small}
\begin{eqnarray*}
E\left[\|Y^{(N)}(t)-y(t)\|^2\right] \leq Ch^{{\min\{2(1-\alpha_m),1\}}}.
\end{eqnarray*}
\end{small}
The proof is completed.
\ \ $\Box$

\section{The fast EM method for multi-term Riemann-Liouville SFDE}
\subsection{The fast EM method}
 To construct the fast EM method, it is necessary to introduce the sum-of-exponentials (SOEs) approximation.
\begin{lemma}\cite{Jiang-2017}\label{SOEs approxiamtion}
For a given $\alpha\in(0,1)$, let $\epsilon$ denote tolerance error, $\delta$ denote cut-off time restriction and $T$ denote final time, there are a positive integer $N_{exp}^{(\alpha)}$ and positive constants $\omega_{j}^{(\alpha)}$ and $s_{j}^{(\alpha)}$, $j=1,2,\cdots,N_{exp}^{(\alpha)}$ such that for any $t\in[\delta,T]$
\begin{small}
\begin{eqnarray*}
\left|t^{-\alpha}-\sum_{j=1}^{N_{exp}^{(\alpha)}}\omega_{j}^{(\alpha)}e^{-s_{j}^{(\alpha)}t}\right|\leq\epsilon,
\end{eqnarray*}
\end{small}
where $N_{exp}^{(\alpha)} = \mathscr{O}((\log\epsilon^{-1})(\log\log\epsilon^{-1}+\log(T\delta^{-1}))+(\log\delta^{-1})(\log\log\epsilon^{-1}+\log\delta^{-1}))$.
\end{lemma}

Considering the modified method (\ref{the discrete EM scheme}) at $t_{n+1}$ and applying Lemma \ref{SOEs approxiamtion} to approximate $(t-\tau_N(s))^{-\alpha_i}$ in the integral from $0$ to $t_n$, then equation (\ref{the discrete EM scheme}) can be reformed as follows:
\begin{small}
\begin{eqnarray}\label{the discrete fast EM scheme}
&&X^{(N)}(t_{n+1}) \nonumber\\
&=& y_0- \sum_{i=1}^{m}\frac{1}{\Gamma(1-\alpha_i)}\left[\int_{0}^{t_n}\frac{X^{(N)}(\tau_N(s))}{(t_{n+1}-\tau_N(s))^{\alpha_i}}ds+\int_{t_n}^{t_{n+1}}\frac{X^{(N)}(\tau_N(s))}{(t-\tau_N(s))^{\alpha_i}}ds\right]\nonumber\\
&&+\int_{0}^{t}f(\tau_N(s),X^{(N)}(\tau_N(s)))ds + \int_{0}^{t}g(\tau_N(s),X^{(N)}(\tau_N(s)))dW(s)\nonumber\\
&\approx&y_0-\sum_{i=1}^{m}\frac{1}{\Gamma(1-\alpha_i)}\left[\sum_{j=1}^{N_{exp}^{(\alpha_i)}}\omega_{j}^{(\alpha_i)}U_{j}^{(\alpha_i)}(t_{n+1})+h^{1-\alpha_i}X^{(N)}(t_n)\right]\nonumber\\
&&+\sum_{j=0}^{k}f(t_j,X^{(N)}(t_j))h + \sum_{j=0}^{k}g(t_j,X^{(N)}(t_j))\triangle W_j,
\end{eqnarray}
\end{small}
where $\triangle W_j = W(t_{j+1})-W(t_j)$, $n=1,2,\cdots,N-1$ (if $n=0$, the modified EM method (\ref{the discrete EM scheme}) is directly applied) and
\begin{small}
\begin{eqnarray*}
U_{j}^{(\alpha_i)}(t_{n+1}) &=& \int_{0}^{t_n}e^{-s_{j}^{(\alpha_i)}(t_{n+1}-\tau_N(s))}X^{(N)}(\tau_N(s))ds.
\end{eqnarray*}
\end{small}
Especially, the property of exponential functions contributes the following key recurrence relations which will be used to compute $U_{j}^{(\alpha_i)}(t_{n+1})$, $n=1,2,\cdots,N-1$,
\begin{small}
\begin{eqnarray*}
U_{j}^{(\alpha_i)}(t_{n+1}) &=& e^{-hs_{j}^{(\alpha_i)}}U_{j}^{(\alpha_i)}(t_n) + e^{-2hs_{j}^{(\alpha_i)}}X^{(N)}(t_{n-1})h.
\end{eqnarray*}
\end{small}
Similar to (\ref{EM-scheme}), for $t\in(\frac{nT}{N},\frac{(n+1)T}{N}], n=0,1,2,\cdots,N-1$, we introduce the continuous-time version
\begin{small}
\begin{eqnarray}\label{fast EM scheme}
X^{(N)}(t)
 &=& y_0-\sum_{i=1}^{m}\frac{1}{\Gamma(1-\alpha_i)}\int_{0}^{\tau_N(t)}\sum_{j=1}^{N_{exp}^{(\alpha_i)}}\omega_{j}^{(\alpha_i)}e^{-s_{j}^{(\alpha_i)}(t-\tau_N(s))}X^{(N)}(\tau_N(s))ds
\nonumber\\
&&-\sum_{i=1}^{m}\frac{1}{\Gamma(1-\alpha_i)}\int_{\tau_N(t)}^{t}(t-\tau_N(s))^{-\alpha_i}X^{(N)}(\tau_N(s))ds\nonumber\\
&&+\int_{0}^{t}f(\tau_N(s),X^{(N)}(\tau_N(s)))ds + \int_{0}^{t}g(\tau_N(s),X^{(N)}(\tau_N(s)))dW(s).
\end{eqnarray}
\end{small}

\begin{remark}
According to Lemma \ref{SOEs approxiamtion} and the above key recurrence relations, we notice that the fast EM method needs less computational cost than the modified EM method, with the modified EM method and the fast EM method need $\mathscr{O}(N^2)$ and $\mathscr{O}(NN_{exp}^{(\alpha_i)})$ computational cost, respectively and $N_{exp}^{(\alpha_i)}$ is much less than $N$ when the step size is quite small. About more details, we can refer to the article \cite{Jiang-2017}.
\end{remark}

\subsection{Strong convergence}
\begin{theorem}\label{fast_bounded}
If Assumption \ref{assumption_3} holds, then for all $N\in N^{*}$, there exists a positive constant $C$ independent of $N$ such that
\begin{small}
\begin{eqnarray*}
E\left[\|X^{(N)}(t)\|^2\right] \leq C \quad {\rm{and}} \quad E\left[\|X^{(N)}(t_n)\|^2\right] \leq C, \quad \forall t, t_n \in [0,T].
\end{eqnarray*}
\end{small}
\end{theorem}
\noindent\emph{Proof.}
We first notice that it follows from equation (\ref{fast EM scheme}) with Lemma \ref{inequality} that
\begin{small}
\begin{eqnarray*}
&&\frac{1}{5}E\left[\|X^{(N)}(t)\|^2\right]\\
&\leq& E\left[\|y_0\|^2\right] + +m\sum_{i=1}^{m}E\left[\|\frac{1}{\Gamma(1-\alpha_i)}\int_{0}^{t}(t-\tau_N(s))^{-\alpha_i}X^{(N)}(\tau_N(s))ds\|^2\right]\\
&&+m\sum_{i=1}^{m}E\left[|\frac{1}{\Gamma(1-\alpha_i)}\int_{0}^{\tau_N(t)}(\sum_{j=1}^{N_{exp}^{(\alpha_i)}}\omega_{j}^{(\alpha_i)}e^{-s_{j}^{(\alpha_i)}(t-\tau_N(s))}-(t-\tau_N(s))^{-\alpha_i})X^{(N)}(\tau_N(s))ds\|^2\right]\\
&&+E\left[\|\int_{0}^{t}f(\tau_N(s),X^{(N)}(\tau_N(s)))ds\|^2\right] + E\left[\|\int_{0}^{t}g(\tau_N(s),X^{(N)}(\tau_N(s)))dW(s)\|^2\right].
\end{eqnarray*}
\end{small}
Similar to the proof of Theorem \ref{bound}, using H$\ddot{o}$lder inequality, It$\hat{o}$ isometry, Lemma \ref{SOEs approxiamtion} and Assumption \ref{assumption_3} yields
\begin{small}
\begin{eqnarray*}
\frac{1}{5}E\left[\|X^{(N)}(t)\|^2\right] &\leq& E\left[\|y_0\|^2\right] + \sum_{i=1}^{m}\frac{m\tau_N(t)}{\Gamma^2(1-\alpha_i)}\int_{0}^{\tau_N(t)}\epsilon^2E\left[\|X^{(N)}(\tau_N(s))\|^2\right]ds\\
&&+\sum_{i=1}^{m}\frac{mt^{1-\alpha_i}}{\Gamma(1-\alpha_i)\Gamma(2-\alpha_i)}\int_{0}^{t}(t-s)^{-\alpha_i}E\left[\|X^{(N)}(\tau_N(s))\|^2\right]ds\\
&&+(2t+2)L_3^2\int_{0}^{t}1+E\left[\|X^{(N)}(\tau_N(s))\|^2\right]ds\\
&\leq&(E\left[\|y_0\|^2\right]+C) + C(\epsilon^2+1)\int_{0}^{t}(t-s)^{-\alpha}E\left[\|X^{(N)}(\tau_N(s))\|^2\right]ds,
\end{eqnarray*}
\end{small}
where the value of $\alpha$ has been defined in the proof of Theorem \ref{bound}. Finally, the Gronwall's inequality completes the proof. At the same time, the proof of $E\left[|X^{(N)}(t_n)|^2\right] \leq C$ can be similarly got. This completes the proof.
\ \ $\Box$

\begin{lemma}\label{t_t_n_2}
If Assumption \ref{assumption_3} holds, then for all $N\in N^{*}$, there exists a positive constant $C$ independent of $N$ such that for all $t\in(t_n,t_{n+1}], t_n\in[0,T]$,
\begin{small}
\begin{eqnarray*}
E\left[\|X^{(N)}(t)-X^{(N)}(t_n)\|^2\right] \leq C\left(h^{{\min\{2(1-\alpha_m),1\}}}+\epsilon^2\right).
\end{eqnarray*}
\end{small}
\end{lemma}
\noindent\emph{Proof.}
For $t\in(t_n,t_{n+1}], t_n\in[0,T]$, the estimate
\begin{small}
\begin{eqnarray}\label{t_t_n_2_inequality}
&&\frac{1}{6}E\left[\|X^{(N)}(t)-X^{(N)}(t_n)\|^2\right]\nonumber\\
&\leq&m\sum_{i=1}^{m}E\left[\|\frac{1}{\Gamma(1-\alpha_i)}\int_{0}^{\tau_N(t)}(\sum_{j=1}^{N_{exp}^{(\alpha_i)}}\omega_{j}^{(\alpha_i)}e^{-s_{j}^{(\alpha_i)}(t-\tau_N(s))}-(t-\tau_N(s))^{-\alpha_i})X^{(N)}(\tau_N(s))ds\|^2\right]\nonumber\\
&&+m\sum_{i=1}^{m}E\left[\|\frac{}{\Gamma(1-\alpha_i)}\int_{0}^{\tau_N(t_n)}(\sum_{j=1}^{N_{exp}^{(\alpha_i)}}\omega_{j}^{(\alpha_i)}e^{-s_{j}^{(\alpha_i)}(t_n-\tau_N(s))}-(t_n-\tau_N(s))^{-\alpha_i})X^{(N)}(\tau_N(s))ds\|^2\right]\nonumber\\
&&+m\sum_{i=1}^{m}E\left[\|\frac{1}{\Gamma(1-\alpha_i)}\int_{0}^{t_n}[(t-\tau_N(s))^{-\alpha_i}-(t_n-\tau_N(s))^{-\alpha_i}]X^{(N)}(\tau_N(s)ds\|^2\right]\nonumber\\
&&+m\sum_{i=1}^{m}E\left[\|\frac{1}{\Gamma(1-\alpha_i)}\int_{t_n}^{t}(t-\tau_N(s))^{-\alpha_i}X^{(N)}(\tau_N(s))ds\|^2\right]\nonumber\\
&&+E\left[\|\int_{t_n}^{t}f(\tau_N(s),X^{(N)}(\tau_N(s)))ds\|^2\right] + E\left[\|\int_{t_n}^{t}g(\tau_N(s),X^{(N)}(\tau_N(s)))dW(s)\|^2\right]\nonumber\\
&=& M_1 + M_2 + M_3 + M_4 + M_5 + M_6
\end{eqnarray}
\end{small}
can be derived by using Lemma \ref{inequality}. Using H$\ddot{o}$lder inequality, Lemma \ref{SOEs approxiamtion} and Theorem \ref{fast_bounded} leads to
\begin{small}
\begin{eqnarray*}
M_1 + M_2 \leq \sum_{i=1}^{m}\frac{m({\tau_N(t)}^2+{\tau_N(t_n)}^2)}{\Gamma^2(1-\alpha_i)}C\epsilon^2\leq C\epsilon^2.
\end{eqnarray*}
\end{small}
Then according to the proofs of Lemma \ref{t_t_n} and Theorem \ref{fast_bounded}, we derive
\begin{small}
\begin{eqnarray*}
M_3 + M_4 + M_5 + M_6 \leq Ch^{{\min\{2(1-\alpha_m),1\}}}.
\end{eqnarray*}
\end{small}
Inserting the two above inequalities into (\ref{t_t_n_2_inequality}) gives
\begin{small}
\begin{eqnarray*}
E\left[\|X^{(N)}(t)-X^{(N)}(t_n)\|^2\right] \leq Ch^{{\min\{2(1-\alpha_m),1\}}}+C\epsilon^2.
\end{eqnarray*}
\end{small}
This completes the proof.
\ \ $\Box$

\begin{theorem}(Strong convergence)\label{fast convergence}
Under Assumptions \ref{assumption_1}-\ref{assumption_3}, there exists a positive constant $C$ independent of $N$ such that for all $t\in[0,T]$,
\begin{small}
\begin{eqnarray*}
E\left[\|X^{(N)}(t)-y(t)\|^2\right] \leq C\left(h^{{\min\{2(1-\alpha_m),\ 1\}}}+\epsilon^2\right).
\end{eqnarray*}
\end{small}
\end{theorem}
\noindent\emph{Proof.}
The desired result can be obtained based on the above lemmas and this theorem can be similarly proved like the Theorem \ref{convergence}.
\ \ $\Box$

\section{Numerical experiments}
In this section, we give two numerical examples to verify the theoretical results of the modified EM method (\ref{EM-scheme}) and the fast EM method (\ref{fast EM scheme}) and demonstrate the performance of the two methods. It can be observed from Theorem \ref{convergence} and Theorem \ref{fast convergence} that both the two proposed meyhods are $\min\{2(1-\alpha_i),\ 1\}_{i=1}^{m}$. Moreover, from the theoretical analysis, the fast EM method has overwhelming advantages over the modified EM method in computational efficiency. In this paper, all of computations are performed by using a MATLAB (R2021b) subroutine on a desktop computer (Dell OptiPlex 7490 AIO) with the Intel(R) Core(TM) i7-11700 CPU @2.50 GHz and 32.0G RAM. In a similar way as \cite{Cao-2015}, the expectation is approximated by sample average. In the following simulations, we define the mean-square errors as
\begin{eqnarray*}
e_n = \max_{1 \leq k \leq n}\left(\frac{1}{5000}\sum_{i=1}^{5000}\|Y^{(n)}(t_k,\omega_i)-Y^{(2n)}(t_k,\omega_i)\|^2\right)^{1/2},
\end{eqnarray*}
where $\omega_i$ represents the $i$th sample path, $Y^{(n)}(t_k,\omega_i)$ and $Y^{(2n)}(t_k,\omega_i)$ denote the $i$th numerical solution and the $i$th numerical solution, respectively.

\begin{example}\label{example_1}
Consider the following two-term SFDEs
\begin{small}
\begin{eqnarray}\label{example_1_equation}
&&y'(t)+\sum_{i=1}^{2}\sideset{^{RL}_0}{^{\alpha_i}_t}{\mathop{\mathrm{D}}}y(t) = sin(t(y(t))) + sin(y(t))\frac{dW_t}{dt}, t\in(0,1],
\end{eqnarray}
\end{small}
with initial value $y(0)=y_0=0.1$.
\end{example}

It can be easily verified that the functions $f$ and $g$ satisfy the assumptions of Theorems \ref{convergence} and \ref{fast convergence}. To check the convergence order of the modified EM method (\ref{EM-scheme}) and the fast EM method (\ref{fast EM scheme}), we use the proposed two methods to compute the errors and convergence orders for $n=128, 256, 512, 1024$ with different combination of $\alpha_1$ and $\alpha_2$. According to the results of these computations listed in Tables 1 and 2, we can conclude that both the convergence orders of the two methods are close to $\min\{1-\alpha_2,\ 0.5\}$, which are consistent with our theoretical results. Meanwhile, from the average CPU times in Tables 1 and 2, it is clear that the CPU time of the fast method is extremely less that that's of the modified EM method especially for the small step sizes.
\begin{table}[!htbp]
 \caption{For $\alpha_1 = 0.6, \alpha_2 = 0.8$, convergence orders of the modified EM method (\ref{EM-scheme}) and the fast EM method (\ref{fast EM scheme}) for Example \ref{example_1}.}
\begin{center}
 \begin{tabular}{lcccccccccc}
  \hline
   \multirow{2}{*}{$n$} &  \multicolumn{3}{c}{modified EM method} & &  \multicolumn{3}{c}{fast EM method}\\
   \cline{2-4}\cline{6-8}     &  error     &  order    & CPU Time& &  error     &  order  & CPU Time  \\
  \hline
   128   &3.543e-3 & &8.66 & &3.543e-3 & &2.15 \\
   256   &3.131e-3 &0.178 &35.01& &3.131e-3 &0.178 &5.61\\
   512   &2.798e-3 &0.162 &143.03& &2.798e-3 &0.162 &15.87\\
   1024  &2.486e-3 &0.171 &559.27& &2.486e-3 &0.171 &51.94\\
  \hline
 \end{tabular}
 \end{center}
\end{table}

\begin{table}[!htbp]
 \caption{For $\alpha_1 = 0.1, \alpha_2 = 0.2$, convergence orders of the modified EM method (\ref{EM-scheme}) and the fast EM method (\ref{fast EM scheme}) for Example \ref{example_1}.}
\begin{center}
 \begin{tabular}{lcccccccccc}
  \hline
   \multirow{2}{*}{$n$} &  \multicolumn{3}{c}{modified EM method} & &  \multicolumn{3}{c}{fast EM method}\\
   \cline{2-4}\cline{6-8}     &  error     &  order    & CPU Time& &  error     &  order  & CPU Time  \\
  \hline
   128   &1.397e-3 & &8.30 & &1.397e-3 & &2.06\\
   256   &9.982e-4 &0.485 &32.91& &9.982e-4 &0.485 &5.40\\
   512   &7.092e-4 &0.493 &131.29& &7.092e-4 &0.493 &16.65\\
   1024  &5.174e-4 &0.455 &525.11& &5.174e-4 &0.455 &51.98\\
  \hline
 \end{tabular}
 \end{center}
\end{table}

In order to better verify the theoretical results of this article, another numerical example of three-term SFDEs is given below.

\begin{example}\label{example_2}
Consider the following three-term nonlinear SFDEs
\begin{small}
\begin{eqnarray}\label{example_1_equation}
\sideset{^{RL}_0}{^{\alpha_1}_t}{\mathop{\mathrm{D}}}y(t)+\sideset{^{RL}_0}{^{\alpha_2}_t}{\mathop{\mathrm{D}}}y(t)+\sideset{^{RL}_0}{^{\alpha_3}_t}{\mathop{\mathrm{D}}}y(t) = sin(t(y(t))) + sin(y(t))\frac{dW_t}{dt},
\end{eqnarray}
\end{small}
for $t\in(0,1]$ and the initial value $y(0)=y_0=0.1$.
\end{example}

Tables 3 and 4 present the errors, convergence orders and CPU times of the two EM methods (\ref{EM-scheme}) and (\ref{fast EM scheme}) with different combinations of $\alpha_1$, $\alpha_2$ and $\alpha_3$. In Tables 3 and 4, we let $\alpha_1 = 0.1, \alpha_2 = 0.5, \alpha_3 = 0.85$ and $\alpha_1 = 0.3, \alpha_2 = 0.35, \alpha_3 = 0.4$, respectively. And the results in Tables 3 and 4 show that the convergence orders of the two methods approach to 0.15 and 0.5, respectively, which are consistent with the theoretical convergence order $\min\{1-\alpha_3,\ 0.5\}$. What's more, we notice that the computational performance of the fast EM method has the overwhelming advantages over the modified EM method from the records of CPU times.

\begin{table}[!htbp]
 \caption{For $\alpha_1 = 0.1, \alpha_2 = 0.5, \alpha_3 = 0.85$, convergence orders of the modified EM method (\ref{EM-scheme}) and the fast EM method (\ref{fast EM scheme}) for Example \ref{example_2}.}
\begin{center}
 \begin{tabular}{lcccccccccc}
  \hline
   \multirow{2}{*}{$n$} &  \multicolumn{3}{c}{modified EM method} & &  \multicolumn{3}{c}{fast EM method}\\
   \cline{2-4}\cline{6-8}     &  error     &  order    & CPU Time& &  error     &  order  & CPU Time  \\
  \hline
   128   &3.192e-3 & &12.18 & &3.192e-3 & &1.29 \\
   256   &2.911e-3 &0.133 &49.06& &2.911e-3 &0.133 &6.16\\
   512   &2.675e-3 &0.122 &195.45& &2.675e-3 &0.122 &17.41\\
   1024  &2.447e-3 &0.129 &784.22& &2.447e-3 &0.129 &55.62\\
   2048  &2.245e-3 &0.124 &3124.02& &2.245e-3 &0.124 &196.71\\
  \hline
 \end{tabular}
 \end{center}
\end{table}

\begin{table}[!htbp]
 \caption{For $\alpha_1 = 0.3, \alpha_2 = 0.35, \alpha_3 = 0.4$, convergence orders of the modified EM method (\ref{EM-scheme}) and the fast EM method (\ref{fast EM scheme}) for Example \ref{example_2}.}
\begin{center}
 \begin{tabular}{lcccccccccc}
  \hline
   \multirow{2}{*}{$n$} &  \multicolumn{3}{c}{modified EM method} & &  \multicolumn{3}{c}{fast EM method}\\
   \cline{2-4}\cline{6-8}     &  error     &  order    & CPU Time& &  error     &  order  & CPU Time  \\
  \hline
   128   &1.366e-3 & &12.69 & &1.366e-3 & &2.61\\
   256   &9.671e-4 &0.498 &32.91& &9.671e-4 &0.498 &6.29\\
   512   &6.909e-4 &0.485 &131.29& &6.909e-4 &0.485 &18.01\\
   1024  &4.917e-4 &0.491 &525.11& &4.917e-4 &0.491 &56.91\\
   2048  &3.457e-4 &0.508 &3127.83& &3.457e-4 &0.508 &198.08\\
  \hline
 \end{tabular}
 \end{center}
\end{table}

\section{Conclusion}
This paper presents a modified EM method to solve multi-term Riemann-Liouville SFDEs and strictly proves the proposed method to be strong convergent with order $\min\{1-\alpha_m,0.5\}$. In view of the large computation cost of the modified EM method, a fast EM method is proposed based on the SOEs approximation to improve the computational efficiency without changing the convergence order. Furthermore, numerical examples verify the theoretical convergence order and demonstrate that the fast EM method has overwhelming computational efficiency compared with the modified EM method.

\section*{Acknowledgments}
This research is supported by the Major Project on New Generation of Artificial Intelligence from MOST of China (Grant No. 2018AAA0101002), and National Natural Science Foundation of China (Grant Nos. 12171466 and 11701502).
\bibliographystyle{abbrv}
\bibliography{references}

\end{document}